 \documentclass[leqno]{amsart}

\usepackage{geometry}
\usepackage{verbatim}
\usepackage[latin1]{inputenc}
\usepackage[T1]{fontenc} 
\usepackage{hyperref}
\usepackage{srcltx}
\usepackage[english]{babel}
\usepackage{mathrsfs}
\usepackage{amscd,amsmath,amssymb,amsthm,amsfonts,epsfig,graphics}

\usepackage{amsmath} 

\usepackage{mathrsfs}
\usepackage{upgreek}

\usepackage{braket}

\usepackage{mathtools} 

\DeclarePairedDelimiter{\abs}{\lvert}{\rvert}
\DeclarePairedDelimiter{\norm}{\lVert}{\rVert}
\DeclarePairedDelimiter{\uple}{\lparen}{\rparen}

\DeclarePairedDelimiter{\pt}{\lparen}{\rparen}

\newcommand*{\sep}{,\,} 

\mathtoolsset{showmanualtags}

\newcommand*{\intcc}[3][]{\mathopen{#1 [} #2 , #3 \mathclose{#1 ]}}

\newcommand*{\eqdef}{\vcentcolon=}

\newcommand{\eu}{\mathrm{e}}

\newcommand{\R}{\mathbb{R}}
\newcommand{\Z}{\mathbb{Z}}
\newcommand{\T}{\mathbb{T}}
\newcommand{\N}{\mathbb{N}}
\newcommand{\X}{\times}
\newcommand{\D}{\partial}

\DeclareMathOperator{\nep}{\eu}

\newcommand*\dif{\mathop{}\!\mathrm{d}}
\newcommand*{\deriv}[3][]{\frac{\dif^{#1}#2}{\dif #3^{#1}}} 

\newcommand{\ve}{\boldsymbol} 
\newcommand*{\fil}[1]{\overline{#1}} 
\newcommand*{\cs}{\mathrm} 

\newtheorem{theorem}{Theorem}[section]
\newtheorem{lemma}[theorem]{Lemma}
\newtheorem{proposition}[theorem]{Proposition}

\newtheorem{definition}[theorem]{Definition}

\usepackage{xcolor}

\newcommand*{\bp}{\ve p}
\newcommand*{\bq}{\ve q}
\newcommand*{\be}{\ve e}

\newcommand*{\bu}{\ve u}
\newcommand*{\bv}{\ve v}
\newcommand*{\bV}{\ve V}
\newcommand*{\bw}{\ve w}
\newcommand*{\bx}{\ve x}

\newcommand*{\Id}{\mathrm{I}}

\begin{document}
\author{Luca Bisconti, Davide Catania} 

\subjclass[2000]{35Q35, 35Q30, 37L30, 76D03, 76F20, 76F65}
 \keywords{Boussinesq equations, Navier-Stokes
  equations, turbulent flows, Large Eddy Simulation (LES),
  deconvolution models, inertial manifold, global attractor.}

\title[On the inertial manifold for the Boussinesq deconvolution
model]{On the existence of an inertial manifold for a deconvolution model of
  the 2D mean Boussinesq equations}

\maketitle

\centerline{\scshape Luca Bisconti } \medskip {\footnotesize
  \centerline{Universit\`a degli Studi di Firenze}
  \centerline{Dipartimento di Matematica e Informatica ``U. Dini''}
  \centerline{Via S. Marta 3, I-50139 Firenze, Italia}
} 

\medskip

\centerline{\scshape Davide Catania} \medskip {\footnotesize
  \centerline{Universit\`a degli Studi di Brescia} \centerline{Sezione
    Matematica (DICATAM)} \centerline{Via Valotti 9, I-25133 Brescia,
    Italia} }

\begin{abstract}
 We show the existence of an inertial manifold
  (i.e. a globally invariant, exponentially attracting,
  finite-dimensional manifold) for the approximate deconvolution model of
  the 2D mean Boussinesq equations.  This model is obtained by
  means of the Van Cittern approximate deconvolution operators, which
  is applied to the 2D filtered Boussinesq equations.
\end{abstract}

\section{Introduction} In what follows we prove the existence of an
inertial manifold for a regularized version of the two-dimensional
viscous Boussinesq system
  \begin{align} 
    & \D_t \bu -\nu \Delta \bu + \nabla \cdot(\bu \otimes \bu) +
      \nabla \pi = \theta \be_2, \label{B-vero-1}\\
    &     \D_t \theta -\kappa \Delta \theta + \bu \cdot \nabla \theta
      = f, \label{B-vero-2}\\
    &     \nabla \cdot \bu = 0, \\
    & (\bu, \theta )\vert_{t=0} = (\bu_0, \theta_0), \label{B-vero-4}
  \end{align}
where $\nu>0$ and $\kappa>0$ are viscosities, $\bu = (u_1, u_2)$ is
the velocity field, $\theta$ may be interpreted physically as a
density, $\pi$ is the cinematic pressure,
and $\be_2:= (0, 1)^T$ where $\{\be_1, \be_2\}$ is the canonical basis
of $\R^2$. Here, we also consider an external action $f=f(\bx)$,
independent of time, forcing the evolution of $\theta$.  

The Boussinesq equations are employed, for instance, as a mathematical
scheme to describe Newtonian fluids whenever density, salinity
concentration, or temperature variations (according to the meaning of
$\theta$) play significant
roles  (see, e.g., \cite{Be-Ca-2013, Be-Ca-2014, Bis} for some recent papers on this subject). 

In performing our analysis we follow closely what has been done in
\cite{HGTiti} where the existence of inertial manifolds has been
proved in the case of two regularizing models for the 2D Navier-Stokes 
equations (NSEs). Differently from reference \cite{HGTiti}, here we
  consider a turbulence model that converges to 2D mean
Navier-Stokes-like equations, namely the 2D mean Boussinesq system, 
when the deconvolution parameter $N$ goes to infinity (see
\eqref{eq:deconv-par}, below). 
Further, we take into account a double filtered system,
  that is a model in which both equations for $\bu$ and $\theta$ 
are somehow regularised by an application of the inverse of the  
Helmholtz operator (see \eqref{eq.w}--\eqref{data} below).

We set $A=\rm \Id-\alpha^2\Delta$, $\fil{\bu}=A^{-1}\bu$, and take
$g(\bx)=\fil{f}(\bx)= A^{-1}f (\bx)$. Filtering the equations
\eqref{B-vero-1}--\eqref{B-vero-4}, we obtain what we call the 
``mean Boussinesq equations'', i.e.
\begin{equation*} 
\left. \begin{array}{l} \D_t
      \overline{\bu} -\nu\Delta \overline{\bu} + \nabla\cdot
      \overline{(\bu \otimes \bu)} + \nabla \overline{\pi} =
      \overline{\theta} \be_2,\\
      \D_t \overline{\theta} -\kappa\Delta \fil{\theta} +
      \nabla \cdot \overline{(\theta \bu)} = \fil{f},\\
      \nabla \cdot \overline{\bu} = 0, \\
      (\overline{\bu}, \overline{\theta} )\vert_{t=0} =
      (\overline{\bu_0}, \overline{\theta_0}),
    \end{array}\right.
\end{equation*}
where $\bu\otimes \bu := (u_1 \bu , u_2 \bu )$, and 
we supply this problem, which is equivalent to
\eqref{B-vero-1}--\eqref{B-vero-4}, 
with periodic boundary conditions (i.e., the torus $\Omega=\T^2$ is the considered domain).

Here, we consider the Approximate Deconvolution Model (ADM),
introduced by Adams and Stolz \cite{Adam-Stolz, St-Ad-1999,
  St-Ad-2001} (see also  \cite{Lew-2009, Be-Lew-2012,  Be-Ca-Lew-2013,
  Bis});  
by following this scheme we approximate the
filtered bi-linear terms as follows:
\begin{equation*}
  \overline{(\bv \otimes \bv)} \sim 
  \overline{(D_N (\overline{\bv}) \otimes D_N (\overline{\bv}))}
  \textrm{ and } 
  \overline{(\varphi \bv)} \sim 
  \overline{(D_N (\overline{\varphi}) D_N (\overline{\bv}))},
\end{equation*}
where $\bv$ and $\varphi$ play the role of the variables $\bu$ and
$\theta$, respectively, and the filtering operator $G_\alpha$ is
defined by the Helmholtz filter (see, e.g., \cite{Be-Ili-Lay-2006,  Be-Lew-2012, 
  Lew-2009, BerBis}, see also \cite{BisCat, BisCat2} for an analogous
case involving an anisotropic horizontal filter), with $\overline{(\,\cdot
  \,)}= G_\alpha(\, \cdot\, )$ and
$G_\alpha:=(I-\alpha^2\Delta)^{-1}$. Here, $D_N$ is the deconvolution
operator, which is constructed using the Van~Cittert algorithm (see,
e.g, \cite{Lew-2009}) and is formally defined by
\begin{equation} \label{eq:deconv-par} 
D_N := \sum_{n=0}^N (I -  G_{\alpha})^n \textrm{ with } N\in\N.
\end{equation}

The ADM we analyze here  is defined by 
\begin{equation*} 
B(\bw, \bw) :=
  \overline{D_N(\bw) \otimes D_N(\bw)}, \,\,\, \mathcal{B}(\rho, \bw)
  := \overline{D_N(\rho) D_N(\bw)},
\end{equation*}
and the considered system  (see \cite{Bis, Be-Lew-2012} for the global 
existence and well-posedness  of the 3D case) is the following
\begin{align}
  &\partial_t \bw - \nu \Delta\bw + \nabla \cdot \fil{D_N \bw \otimes
    D_N \bw} +\nabla q = \rho \be_2\, , \label{eq.w}\\
  &\partial_t \rho - k \Delta\rho + \nabla \cdot \fil{D_N\rho D_N\bw} = g\, , \label{eq.rho} \\
  &\nabla\cdot\bw = 0\,  \label{eq.div}\\
&  \bw(0, x) = \overline{\bu_0}(x),\,\, \textrm{ and }\,\, \rho(0, t) =
  \overline{\theta_0}(x). \label{data}
\end{align}
This system of partial differential equations (PDEs) 
has a dissipative nature and so, energy associated with weak solutions
decreases during time. Consequently, we can presume the existence of a 
global attractor $\mathcal{A}$, i.e. a compact and non-empty subset of the
phase space $\mathcal{H}$ (with metric $d$), of the initial data, 
which is invariant with respect to the action of the semigroup $S(t)$
associated with the system itself (see, e.g., \cite{Temam-88, C-V-85, Flandoli, Zelik}). 
Again,  $\mathcal{A}$ has the property of attracting images, under the action of $S(t)$, of any
bounded subsets $B\subseteq \mathcal{H}$  when time $t$
goes to $+\infty$ (i.e., $\lim_{t\to +\infty}  \mathrm{d}(S(t) B,
\mathcal{A})=0$). Consequently, $\mathcal{A}$ consists of
the full trajectories of the considered system and contains
all the related ``non-trivial'' dynamics (see, e.g., \cite{Zelik}).

 In fact, in the case of system \eqref{eq.w}--\eqref{data}, we can say more because the associated
dynamics is effectively finite dimensional and it can be completely
analyzed by using a suitable system of ordinary differential equations
(ODEs). We prove this fact in the main result of the paper,
i.e. Theorem~\ref{thm.main}, by considering suitably modified equations (see later on).

By definition (see \cite{F-1, F-2}), an inertial
manifold for an evolution equation is a finite-dimensional smooth
submanifold $\mathcal M$ of the phase space
$\mathcal{H}$, which is invariant with respect to the action of the
semigroup $S(t)$, it contains the global
attractor $\mathcal{A}$, and it is also such that any trajectory
starting outside of $\mathcal M$ approaches, exponentially
fast, some other trajectory belonging to $\mathcal M$. If this object
  exists, by restricting
the evolution equation (here, in fact, we think of a single or a system of PDEs) 
to $\mathcal M$, we get a finite system of ODEs that captures the dynamics
on the attractor.  This system of
ODEs is called the ``inertial form of the given evolutionary
system''. As a result, the dynamical properties of the solution of the
evolutionary PDE (or of the system of PDEs), which is an infinite-dimensional system,
can be analyzed by treating a finite-dimensional scheme,
which is the corresponding inertial form.

As remarked in \cite{HGTiti}, the original motivations for the
development of the theory of inertial manifolds were mainly related to
the analysis of  the NSEs. Nevertheless, as far as we know,
the problem about the existence of inertial manifolds for the 2D NSEs
is still open. As a partial answer in this direction, the existence of inertial manifolds has been
proved for some ``$\alpha$-models''; to be precise, for the 2D Bardina
model and for the 2D modified Leray-$\alpha$ 
(see \cite{HGTiti}). \smallskip

\noindent The paper is organized as follows: in Section~\ref{basics} we
introduce preliminary notions and the functional setting. In
Section~\ref{estimates} we provide some a-priori estimates
 that highlight the dissipative nature of the considered equations.
 Section~\ref{InertialMan} is devoted to prove the
existence of absorbing sets in suitable Hilbert spaces, as well as to
show the existence of an inertial manifold for a modification of the approximate deconvolution  system
\eqref{eq.w}--\eqref{data}.

\section{Basic facts and notation} \label{basics} 
We introduce the functional setting, and we recall the definition and 
the main properties of the deconvolution operator $D_N$.

Denote by $x := (x_1, x_2) \in \R^2$ a generic point in
$\R^2$. Given $L \in \R^\star_+ := \{x \in \R : x > 0\}$,
consider $ ]0, L[^2 \subset \R^2$ and ${\mathcal{T}}_2 :=
2\pi\Z^2/L$; $\Omega=\T^2$ is the torus defined by $\T^2 := \big(\R^2/
{\mathcal{T}}_2\big)$.  
For $p, k \in \N$, we consider the standard Lebesgue spaces
$L^p=L^p_{per}=L^p(\T^2)$ and Sobolev spaces $W^{k, p}=W^{k, p}(\T^2)$, with
$H^k:=W^{k, 2}$, in the periodic setting.  We
denote by $\|\cdot\|$ the $L^2(\T^2)$-norm and the associated operator
norms, and we always impose the zero mean condition on the considered
fields. In the sequel, we will use the same notation for scalar and
vector-valued functions, since no ambiguity occurs.  Also, dealing
with divergence-free vector fields, we also define, for a generic
exponent $s \geq 0$, the following spaces
\begin{equation*}
  H_s :=
  \Big\{ \bv : \T^2 \to \R^2\, \colon\,  \bv \in (H^s)^2,\,\, 
  \nabla \cdot \bv = 0,
  \,\,\int_{\T^2}\bv dx = 0\Big\}.
\end{equation*}
When $0\leq s \leq 1$, the condition $\nabla\cdot \bv =0$ must be
understood in weak sense. Let $X$ be a real Banach space with norm
$\|\cdot\|_X$.  We will use the customary Bochner spaces $L^q(0, T;
X)$, with norm denoted by $\|\cdot\|_{L^q(0,T;X)}$.
We denote by $P_\sigma : L^2_{per} \to H:=H^0$ the Helmholtz-Leray 
orthogonal projection operator, and by $\mathbf{A} = -P_\sigma \Delta$ 
the Stokes operator. 
Since we are dealing with periodic-boundary conditions, we have that
\begin{equation} \label{projector}
\mathbf{A}u = -P_{\sigma}\Delta u = -\Delta u,\,\, \textrm{for all $u \in D(\mathbf{A})$}.
\end{equation}

For $v \in H^s$, we can expand the fields as $v(x) = \sum_{k\in
  {\mathcal{T}}_2^\star} \widehat{v}_k e^{ik\cdot x}$, where $k \in
{\mathcal{T}}_2^\star$, and the Fourier coefficients are defined by
$\widehat{v}_k = 1/|\T^2|\int_{{\T}^2} v(x) e^{-ik \cdot x}dx$.  The
magnitude of $k$ is given by $|k|^2 := (k_1)^2 + (k_2)^2$. 
The $H^s$ norms are defined by $\|v\|^2_s:= \sum_{k\in
  {\mathcal{T}}_2^\star} |k|^{2s} |\widehat{v}_k|^2$, where $\|v\|^2_0
:= \|v\|^2$.  The inner products associated to these norms are $(w,
v)_{H^s} := \sum_{k\in {\mathcal{T}}_2^\star}|k|^{2s}\widehat{w}_k
\cdot \overline{\widehat{v}_k}$, where $\overline{\widehat{v}_k}$
denotes the complex conjugate of $\widehat{v}_k$.  To have real valued
vector fields, we impose $\widehat{v}_{-k} = \overline{\widehat{v}_k}$
for any $k \in {\mathcal{T}}_2^\star$ and for any field denoted by
$v$. When $s$ is an integer, $\|v\|^2_s = \|\nabla^s v\|^2$ and also, for general $s
\in \R$, $(H^s)' = H^{-s}$,  (see e.g. \cite{Do-Gib-1995}). All these considerations can be adapted
straightforwardly to the case of the spaces $H_s$; we finally
characterize $H_s \subseteq \big(W^{s,2}(\T^2)\big)^2$ as
follows:
\begin{equation*}
H_s =
  \Big\{ \bv =\sum_{|k|\geq 1} \widehat{\bv}_k e^{ik\cdot x}\colon\,
  \sum_{|k|\geq 1} |k|^{2s} |\widehat{\bv}_k|^2< \infty,\,\, 
  k\cdot \widehat{\bv}_k = 0,
  \,\, \widehat{\bv}_{-k} = \overline{\widehat{\bv}_k}\Big\}.
\end{equation*}
  In particular, we
denote $(H_s)'$ by $H_{-s}$.

We will denote by $C$ generic constants, which may change from line to
line, but which are independent of the diffusion coefficient
$\epsilon$, the deconvolution parameter $N$ and of the solution of the
equations we are considering.

Let us now briefly recall the properties of the Helmholtz filter.
Let $\alpha > 0$ be a given fixed number and, for $\ve \omega \in H_s$, $s\geq
-1$, let us denote by $(\overline{\ve \omega}, \pi) \in H_{s+2} \X H^{s+1}$,
the unique solution of the following Stokes-like problem:
\begin{equation} \label{eq:helmholtz-filter}
  \begin{aligned}
    & \overline{\ve \omega} -\alpha^2\Delta \overline{\ve \omega} +
    \nabla \pi = \ve \omega \textrm{ in } \T^2,\\
    & \nabla \cdot \overline{\ve \omega} = 0 \textrm{ in } \T^2,
  \end{aligned}
\end{equation}
with $\int_{\T^2} \overline{\ve \omega} dx = 0$ and $\int_{\T^2} \pi dx = 0.$
The velocity component of $(\overline{\ve \omega}, \pi)$ is denoted also by
$\overline{\ve \omega} = G_\alpha(\ve \omega)$ and $A_1:=G_\alpha^{-1}$.  
Consider an element $\ve \omega \in H_s$ and take its expansion in terms of
Fourier series as $\ve \omega = \sum_{k\in {\mathcal{T}}_2^\star}\widehat{\ve \omega}_k
e^{ik\cdot x}$, so that inserting this expression in
\eqref{eq:helmholtz-filter} and looking for $(\overline{\omega},
\pi)$, we get 
\begin{equation} \label{eq:pure-utility} \overline{\ve \omega} = \sum_{k\in
    {\mathcal{T}}_2^\star} \frac{1}{1+\alpha^2|k|^2} \widehat{\ve  w}_k
  e^{ik\cdot x} = G_\alpha(\ve \omega), \textrm{ and } \pi=0.
\end{equation}

For a scalar function $\zeta$ we denote by $\overline{\zeta}$ the
solution of the pure Helmholtz problem
\begin{equation} \label{eq:pure-helmholtz} 
-\alpha^2\Delta
  \overline{\zeta} + \overline{\zeta} = \zeta\, \textrm{ in } \T^2,
\end{equation}
where $A_2\overline{\zeta}:= -\alpha^2\Delta \overline{\zeta} +
\overline{\zeta}$. Further, taking $\zeta\in H^s$ the expression of
$\overline{\zeta}$ in terms of Fourier series can be retrieved,
formally, by \eqref{eq:pure-helmholtz} substituting $\zeta$ in place of
$w$. 

In what follows, in order to keep the notation compact, we use the
same symbol $A$ for the operators $A_1$ and $A_2$, distinguishing the
two situations only when it is required by the context.  According to
the above facts, the deconvolution operator $D_N$ 
is actually given by $D_N = \sum_{n=0}^N (I - A^{-1})^n$, $N\in\N$,
with $A$ defined by \eqref{eq:helmholtz-filter}, when it is acting on
the elements of $H_s$ and, by relation \eqref{eq:pure-helmholtz}, in the case
of the scalar functions in $H^s$.

Let us recall that, 
 in the filtered equations \eqref{eq.w}--\eqref{eq.rho}, the symbol ``\hspace{0.05
  cm}$ \overline{\empty{{}^{{}^{\,\,\,\,}}}}$\hspace{0.05 cm}'' (i.e
`` $ \overline{\,\cdot
  \,} $ '') denotes the pure Helmholtz filter, applied component-wise to the
various vector and tensor fields. Referring to the right-hand side of
equation \eqref{eq.w}, since $\be_2$ is a constant
vector, then we have that $G_\alpha (\theta \be_2) = \overline{\theta
  \be_2} = \overline{\theta}\be_2 = G_\alpha (\theta)\be_2$ and $A
(\overline{\theta \be_2}) = A(\overline{\theta})\be_2 $ (where the meaning
of $A$ is understood in the sense stated above). Also, for brevity, in
the sequel we omit the explicit dependence of $G_\alpha$ on $\alpha$
(we write $G$ in place of $G_\alpha$) and we will always use the
notation $G = A^{-1} = (I - \alpha^2\Delta)^{-1}$.

The deconvolution operator $D_N$ is constructed thanks to the
Van~Cittert algorithm; the reader can find an exhaustive description and
analysis of the Van Cittert algorithm and its variants in
\cite{Lew-2009}. Here, we only report the properties needed to
describe the considered model.  Let $\omega\in H_s$ (or $\omega \in
H^s$), starting from the expression \eqref{eq:pure-utility}, we can
write the deconvolution operator in terms of Fourier series as follows
\begin{equation} \label{eq:basics-DN}
  \begin{aligned}  &\widehat{D}_N(k) = \sum_{n=0}^N \left(\frac{\alpha^2
        |k|^2}{1+\alpha^2|k|^2}\right)^n = (1 + \alpha^2 |k|^2
    )\varrho_{N,k} \\
    \varrho_{N,k} =& 1-\left( \frac{ \alpha^2|k|^2}{1 +
        \alpha^2|k|^2}\right)^{N+1} \textrm{ and }\,\, D_N(\omega) =
    \sum_{k\in {\mathcal{T}}_2^\star} \widehat{D}_N (k)
    \widehat{\omega}_k e^{ik\cdot x}.
  \end{aligned}
\end{equation}
The properties verified by $\widehat{D}_N$ are summarized in
\cite[Lemma~2.1]{Be-Lew-2012} and \cite[Lemma~2.2]{Be-Lew-2012}
and, even if we do not use them directly,
 we recall here  the main points for the reader's convenience:
\begin{lemma} \label{lem:utility} For each fixed $k\in
  {\mathcal{T}}_2$,
  \begin{equation}
    \widehat{D}_N(k) \to 1+ \alpha^2|k|^2 = \widehat{A}_k, \textrm{ as } 
    N \to +\infty, \textrm{   even if not uniformly in }  k. 
  \end{equation}
  Further, for each $N\in \N$ the operator $D_N \colon H_s \to H_s$ is
  self-adjoint, it commutes with differentiation, and the following
  properties hold true:
  \begin{align} \allowdisplaybreaks & 1 \leq \widehat{D}_N (k) \leq N
    + 1,\,\, \forall k\in {\mathcal{T}}_2^\star,
    \label{eq:dn-1} \\
    \displaybreak[0] &\widehat{D}_N (k) \cong (N + 1) \frac{1 +
      \alpha^2|k|^2}{\alpha^2 |k|^2}
    \textrm{ for large } |k|, \label{eq:estimate-norm-DN}\\
    &\widehat{D}_N (k) \leq 1 + \alpha^2 |k|^2= \widehat{A}_k, \forall
    k \in {\mathcal{T}}_2^\star,
    \alpha > 0, \\
    & \textrm{the map } \omega \mapsto D_N (\omega) \textrm{is an
      isomorphism s.t. } \|D_N\|_{H_s} = O(N + 1), \forall s \geq 0,\\
    & \underset{N\to +\infty}{\lim} D_N (\omega) = A\omega \textrm{ in
    } H_s\,\, \forall s \in \R \textrm{ and } \omega \in H_{s+2}.
  \end{align}
\end{lemma}
This lemma can be readily extended to the case of spaces $H^s$.

\section{Preliminary estimates} \label{estimates}
In what follows we provide a priori estimates needed to study the
long-time dynamics of the solution of
system  \eqref{eq.w}--\eqref{data}. In particular, we are required to justify the existence
of absorbing balls for the dynamical system induced by equation  \eqref{eq.w}--\eqref{data},
in various spaces of functions.  The controls given in the sequel are formal, but it is
 possible to prove them rigorously, by using and adapted Galerkin
 approximation procedure that now we briefly illustrate.

\subsection{Galerkin scheme}
We give some details about the construction 
  of the approximate solutions, which is however classical 
(for more details see, e.g., \cite{Guo, Lions, Be-Lew-2012, Bis}).  Let be
  given $m \in \N\backslash\{0\}$ and define
  \begin{align*}
    &V^m = \Big\{ \vartheta \in H^1 \, \colon \, \int_{\T^3} \vartheta (x) e^{-ik\cdot
      x}dx = 0,\,\, \forall k
    \textrm{ with } |k| > m \Big\},\\
    &\mathbf{V}_m = \Big\{ \mathbf{w} \in H_1 \, \colon \, \int_{\T^3}
    \mathbf{w}(x) e^{-ik\cdot x}dx = 0,\,\, \forall k \textrm{ with }
    |k| > m \Big\},
  \end{align*}
  and let $\{E_j\}_{j=1,\ldots, d_m}$ and
  $\{\mathbf{E}_j\}_{j=1,\ldots, \delta_m}$ be orthogonal bases of
  $V^m$ and $\mathbf{V}_m$ respectively.  Without loss of generality,
  we can assume that the $E_j$'s are eigen-functions of the operator
  $I -\alpha^2\Delta$ introduced in \eqref{eq:pure-helmholtz} as well
  as the $\mathbf{E}_j$'s are eigen-functions of the Stokes-like
  operator associated to \eqref{eq:helmholtz-filter}.  Further, the
  $E_j$'s and $\mathbf{E}_j$'s are supposed to have unitary norm.  We
  denote by $P_m$ the orthogonal projection from $H^1$ onto $V^m$ and,
  similarly, by $\mathbf{P}_m$ the the orthogonal projection from
  $H_1$ onto $\mathbf{V}_m$.
  For every positive integer $m$, we look for an approximate solution
  of equations \eqref{eq.w}--\eqref{data} of the form
  \begin{equation*}
    \rho_m (t, x) = \sum^{d_m}_{j=1}\rho_{m, j}(t) E_j (x)\, \textrm{ and }\, 
  {\ve  w}_m(t, x) =  \sum^{\delta_m}_{j=1}{\ve w}_{m, j}(t) \mathbf{E}_j(x).
  \end{equation*}
  Thanks to the Cauchy-Lipschitz Theorem, we can prove the existence
  of a unique $C^1$ maximal solution $\big({\ve w}_m(t), \rho_m(t)\big)\in
  \mathbf{V}_m\X V^m$ for all $t \in [0, T_m)$ where $T_m>0$ is the
  maximal existence time, to the system
  \begin{align}
    &\begin{aligned} \label{weak-form-m-1} \int_{\mathbb{T}^2} \D_t
    {\ve  w}_m \cdot \ve v - \int_{\mathbb{T}^2} G\big( D_N({\ve w_m}) \otimes
    &  D_N ({\ve w}_m)\big) :  \nabla \ve v \\
      & + \nu \int_{\T^2} \nabla {\ve w}_m : \nabla \ve v =
      \int_{\mathbb{T}^2}\rho_m \mathbf{e}_2 \cdot \ve v,
    \end{aligned}\\[2 mm]
    & \int_{\mathbb{T}^2} \D_t \rho_m \cdot h - \int_{\mathbb{T}^2}
    G\big( D_N(\rho_m) D_N ({\ve w}_m)\big) \cdot \nabla h + k
    \int_{\mathbb{T}^2} \nabla \rho_m \cdot \nabla h
    =   \int_{\mathbb{T}^2}g \cdot h, \label{weak-form-m-2}
  \end{align}
  for all $(\ve v, h)\in L^2(0, T_m; \mathbf{V}_m)\X L^2(0, T_m; V^m)$. 
Actually $T_m = T$ (see \cite{Bis}), and this  concludes the construction 
of the approximate solutions $\{({\ve w}_m, \rho_m)\}_{m\in \N}$.

In what follows, to keep the notation concise, se set $\ve w= \ve w_m$
and $\rho = \rho_m$. Also, we use the same notation $\|\,\cdot\,\|$
for both $\|\, \cdot\|_{H^0}$ and $\|\, \cdot\|_{H_0}$.

 \subsection{Estimates for $\norm{A^{1/2} D_N^{1/2}\rho(t)}$ and $\norm{A^{1/2}
      D_N^{1/2}\bw(t)}$}
We set $\Lambda \eqdef (-\Delta)^{1/2}$ and note that $\norm{\Lambda
  \ve u} = \norm{\nabla\ve u}$.
We test the Equation~\eqref{eq.rho} against $A D_N\rho$ (this means
that we take the scalar product in $\cs L^2(\Omega)$), with $\Omega=\T^2$,
 and obtain, using
the Cauchy--Schwarz inequality,
\begin{align*}
  \frac{1}{2}\deriv{}{t} \norm{A^{1/2}D_N^{1/2}\rho}^2 + k\norm{\nabla
    A^{1/2}D_N^{1/2}\rho}^2 & \leq \abs{\uple{\Lambda A^{1/2}
      D_N^{1/2}\rho \sep \Lambda^{-1} A^{1/2} D_N^{1/2}g}}
  \\
  & \leq \frac{k}{2} \norm{\nabla A^{1/2}D_N^{1/2}\rho}^2 +
  \frac{1}{2k} \norm{\Lambda^{-1} A^{1/2}D_N^{1/2} g}^2\, .
\end{align*}
Setting $y(t)=\norm{A^{1/2}D_N^{1/2}\rho(t)}^2$ and using the
Poincar\'e inequality, we thus obtain
\begin{align*}
  y'(t)+k\lambda_1 y(t) \leq \frac{1}{k} \norm{\Lambda^{-1} A^{1/2}
    D_N^{1/2} g}^2\, .
\end{align*}
Applying the Gronwall lemma, we deduce
\begin{align*}
  y(t)\leq \frac{1}{k^2\lambda_1} \norm{\Lambda^{-1} A^{1/2} D_N^{1/2}
    g}^2 + \pt*{y(0)-\frac{1}{k^2\lambda_1} \norm{\Lambda^{-1} A^{1/2}
      D_N^{1/2} g}^2}\nep^{-k\lambda_1 t} =: R_1(t)^2\, ,
\end{align*}
so that
\begin{align*}
  &\limsup_{t\to+\infty} y(t) \leq \lim_{t\to +\infty}
  R_1(t)^2=\frac{1}{k^2\lambda_1} \norm{\Lambda^{-1} A^{1/2} D_N^{1/2}
    g}^2 =: \frac{r_1^2}{2}\, .
\end{align*}
This means that $y(t)$ and $R_1(t)^2$ enter a ball of radius $r_1^2$
after long enough time.

Then, we test the Equation~\eqref{eq.w} against $A D_N\bw$, obtaining
\begin{align*}
  \frac{1}{2}\deriv{}{t} \norm{A^{1/2}D_N^{1/2}\bw}^2+\nu\norm{\nabla
    A^{1/2}D_N^{1/2}\bw}^2 & \leq \norm{A^{1/2} D_N^{1/2}\bw} \cdot
  \norm{A^{1/2} D_N^{1/2}\rho}
  \\
  & \leq \norm{\nabla A^{1/2} D_N^{1/2}\bw} \cdot
  \lambda_1^{-1/2}\norm{A^{1/2} D_N^{1/2}\rho}
  \\
  & \leq \frac{\nu}{2} \norm{\nabla A^{1/2} D_N^{1/2}\bw}^2 +
  \frac{1}{2\nu\lambda_1} \norm{A^{1/2} D_N^{1/2}\rho}^2\, .
\end{align*}
Setting $z(t)=\norm{A^{1/2}D_N^{1/2}\bw(t)}^2$ and using the
Poincar\'e inequality, we get
\begin{align*}
  z'(t)+\nu\lambda_1 z(t) \leq \frac{1}{\nu\lambda_1}
  \norm{A^{1/2}D_N^{1/2}\rho(t)}^2 \leq
  \frac{r_1^2+2y(0)\nep^{-k\lambda_1 t}}{2\nu\lambda_1}\, .
\end{align*}
If we multiply both members by $\nep^{\nu\lambda_1 t}$, integrate in
time on $\intcc{0}{t}$ and then divide by $\nep^{\nu\lambda_1 t}$, we
deduce
\begin{align*}
  z(t) \leq z(0)\nep^{-\nu\lambda_1 t} +
  \frac{r_1^2}{2\nu^2\lambda_1^2} (1-\nep^{-\nu\lambda_1 t})
  +\frac{y(0)}{\nu\lambda_1}\nep^{-\nu\lambda_1 t}\int_0^t
  \nep^{(\nu-k)\lambda_1 s}\dif s \, .
\end{align*}
We need to distinguish two cases: $\nu= k$ and $\nu\neq k$. If $\nu=
k$, we have
\begin{align*}
  z(t) \leq \frac{r_1^2}{2\nu^2\lambda_1^2} + \pt*{z(0)-
    \frac{r_1^2}{2\nu^2\lambda_1^2} +
    \frac{y(0)}{\nu\lambda_1}t}\nep^{-\nu\lambda_1 t}\, .
\end{align*}
If $\nu\neq k$, then
\begin{align*}
  z(t) \leq \frac{r_1^2}{2\nu^2\lambda_1^2} + \pt*{z(0)-
    \frac{r_1^2}{2\nu^2\lambda_1^2}}\nep^{-\nu\lambda_1 t}+
  \frac{y(0)}{\nu(\nu-k)\lambda_1^2} \pt*{\nep^{-k\lambda_1
      t}-\nep^{-\nu\lambda_1 t}}\, .
\end{align*}
In any case, we obtain
\begin{align*}
  & z(t) \leq R_2(t)^2 := \begin{cases}
    \dfrac{r_1^2}{2\nu^2\lambda_1^2} + \pt*{z(0)- \dfrac{r_1^2}{2\nu^2\lambda_1^2} + \dfrac{y(0)}{\nu\lambda_1}t}\nep^{-\nu\lambda_1 t} & \quad \text{if } \nu=k\, , \\[1em]
    \dfrac{r_1^2}{2\nu^2\lambda_1^2} + \pt*{z(0)-
      \dfrac{r_1^2}{2\nu^2\lambda_1^2}}\nep^{-\nu\lambda_1 t}+
    \dfrac{y(0)}{\nu(\nu-k)\lambda_1^2} \pt*{\nep^{-k\lambda_1
        t}-\nep^{-\nu\lambda_1 t}} & \quad \text{if } \nu\neq k\, ,
  \end{cases} \\
  & \limsup_{t\to +\infty} z(t) \leq \lim_{t\to +\infty} R_2(t)^2 =
  \frac{r_1^2}{2\nu^2\lambda_1^2} =: \frac{r_2^2}{2}\, .
\end{align*}
This means that $z(t)$ and $R_2(t)^2$ enter a ball of radius $r_2^2$
after long enough time. If we set $r:=\max\set{r_1, r_2}$, we deduce
that there exists $t_r>0$ such that
\begin{equation} \label{eq:time-r}
  \norm{A^{1/2}D_N^{1/2}\rho(t)} \leq R_1(t) \leq r, \,\, \textrm{ and
  }\,\,   \norm{A^{1/2}D_N^{1/2}\bw(t)} \leq R_2(t) \leq r
\end{equation}
for every $t\geq t_r$.

\subsection{Estimates for $\norm{A D_N^{1/2}\rho(t)}$ and $\norm{A
		D_N^{1/2}\bw(t)}$}
In what follows we look for higher order estimates. This time, we start by testing
the Equation~\eqref{eq.w} by $A^2 D_N\bw$ and noticing that the
nonlinear term cancels out (this holds true just in the 2D periodic
case), so that
\begin{align*}
  &\frac{1}{2} \deriv{}{t} \norm{A D_N^{1/2}\bw}^2 + 
\nu\norm{\nabla A D_N^{1/2}\bw}^2 \leq \abs{\uple{A^{1/2}D_N^{1/2}
\rho \sep A^{3/2} D_N^{1/2}\bw}} \\
  & \qquad \leq \frac{\nu}{2} \norm{A^{3/2}D_N^{1/2}\bw}^2 +
  \frac{1}{2\nu} \norm{A^{1/2}D_N^{1/2}\rho}^2\, ,
\end{align*}
and hence
\begin{align*}
  \deriv{}{t} \norm{A D_N^{1/2}\bw}^2 + \nu\norm{\nabla A
    D_N^{1/2}\bw}^2 \leq \frac{R_1^2}{\nu}\, .
\end{align*}
Setting $Z(t)=\norm{A D_N^{1/2}\bw(t)}^2$ and using the Poincar\'e
inequality, we have $Z'(t)+\nu\lambda_1 Z(t)\leq
R_1(t)^2/\nu$. Multiplying by $\nep^{\nu\lambda_1 t}$ and integrating
on $\intcc{t_r}{t}$, with $t\geq t_r$, where $R_1(t)\leq r$, we obtain
\begin{gather*}
  Z(t)  \leq Z(t_r) \nep^{-\nu\lambda_1 (t-t_r)} +
  \frac{r^2}{\nu^2\lambda_1} 
\pt*{1- \nep^{-\nu\lambda_1 (t-t_r)}} \\
  \limsup_{t\to+\infty} Z(t) \leq \frac{r^2}{\nu^2\lambda_1} =:
  \frac{s_2^2}{2}\, .
\end{gather*}

Now, we test the Equation~\eqref{eq.rho} against $A^2 D_N\rho$,
obtaining (here the nonlinear term does not disappear):
\begin{align*}
  \frac{1}{2} \deriv{}{t} \norm{A D_N^{1/2} \rho}^2 + k \norm{A^{3/2}
    D_N^{1/2} \rho}^2 \leq \abs{\uple{A^{1/2} D_N^{1/2} g\sep A^{3/2}
      D_N^{1/2} \rho}} + \abs{\uple{D_N \bw \cdot \nabla D_N \rho \sep
      A D_N \rho}}\, .
\end{align*}
Then we observe that
\begin{align*}
  \abs{\uple{A^{1/2} D_N^{1/2} g\sep A^{3/2} D_N^{1/2} \rho}} \leq
  \frac{k}{4} \norm{A^{3/2} D_N^{1/2}\rho}^2 + \frac{1}{k}
  \norm{A^{1/2} D_N^{1/2} g}^2
\end{align*}
and
\begin{align*}
  & \abs{\uple{D_N \bw \cdot \nabla D_N \rho \sep A D_N \rho}} = 
\abs{\uple{D_N \bw \cdot \nabla A D_N \rho \sep D_N \rho}} \leq
 \norm{D_N \bw}_{\cs L^4} \norm{\nabla A D_N \rho} \, \norm{D_N\rho}_{\cs L^4} \\
  & \qquad \leq c_4 \norm{D_N\bw}^{1/2} \norm{\nabla D_N \bw}^{1/2}
  \norm{D_N\rho}^{1/2} \norm{\nabla D_N \rho}^{1/2} \norm{\nabla A
    D_N\rho}\, .
\end{align*}
Here, we are using the Ladyzhenskaja inequality $\norm{u}_{\cs L^4}
\leq c_4 \norm{u}_{\cs L^2}^{1/2} \norm{\nabla u}_{\cs
  L^2}^{1/2} $ for a scalar function $u$ defined on a torus
$\T^2$. Note that, if $\bu=(u_1,u_2)$ is a vector field, still defined
on $\T^2$,
and we set $\norm{\bu}_{\cs L^p}^p = \norm{u_1}_{\cs L^p}^p +
\norm{u_2}_{\cs L^p}^p $, we have also $\norm{\bu}_{\cs L^4} \leq
c_4 \norm{\bu}_{\cs L^2}^{1/2} \norm{\nabla \bu}_{\cs L^2}^{1/2}
$ (with the same constant $c_4>0$).
Using
\begin{align*}
  \norm{D_N \bw} \leq \norm{A^{1/2} D_N^{1/2} \bw},\,\, \textrm{ and }\,\,
  \norm{\nabla D_N \bw} \leq \frac{(N+1)^{1/2}}{\alpha}\norm{A^{1/2}
    D_N^{1/2} \bw}
\end{align*}
and similar inequalities, we deduce
\begin{align*}
  & \abs{\uple{D_N \bw \cdot \nabla D_N \rho \sep A D_N \rho}} 
\leq c_4^2 \frac{N+1}{\alpha^2} \norm{A^{1/2}D_N^{1/2}\bw}\, 
\norm{A^{1/2}D_N^{1/2}\rho}\, \norm{A^{3/2} D_N^{1/2} \rho} \\
  & \qquad \leq \frac{k}{4}\norm{A^{3/2}D_N^{1/2}\rho}^2 +
  \frac{c_4^4(N+1)^2}{k\alpha^4} \norm{A^{1/2}D_N^{1/2}\bw}^2
  \norm{A^{1/2}D_N^{1/2}\rho}^2 \, .
\end{align*}
Substituting above and multiplying by $2$, we obtain
\begin{align*}
  \deriv{}{t} \norm{A D_N^{1/2} \rho}^2  + k \norm{A^{3/2} D_N^{1/2} \rho}^2 & \leq \frac{2}{k} 
  \norm{A^{1/2} D_N^{1/2} g}^2 + \frac{2c_4^4(N+1)^2}{k\alpha^4} \norm{A^{1/2}D_N^{1/2}\bw}^2 
  \norm{A^{1/2}D_N^{1/2}\rho}^2 \\
  & \leq \frac{2}{k} \norm{A^{1/2} D_N^{1/2} g}^2 + \frac{2c_4^4(N+1)^2}{k\alpha^4} R_1^2 R_2^2 \\
  & \leq \frac{2}{k} \norm{A^{1/2} D_N^{1/2} g}^2 +
  \frac{2c_4^4(N+1)^2}{k\alpha^4} r^4 =: \beta \qquad \text{for }
  t\geq t_r\, .
\end{align*}
Setting $Y(t)=\norm{A D_N^{1/2} \rho}^2$ and applying the Poincar\'e
inequality, we have $Y'(t)+k\lambda_1 Y(t) \leq \beta$ (when $t\geq
t_r$). Proceeding as before, we deduce
\begin{align*}
  & Y(t) \leq  Y(t_r) \nep^{-k\lambda_1 (t-t_r)} + \frac{\beta}{k\lambda_1} 
  \pt*{1- \nep^{-k\lambda_1 (t-t_r)}}\, , \\
  & \limsup_{t\to +\infty} Y(t) \leq \frac{\beta}{k\lambda_1} =:
  \frac{s_1^2}{2}\, .
\end{align*}
If we set $s:=\max\set{s_1, s_2}$, we conclude that there exists
$t_s\geq t_r >0$ such that
\begin{equation*} 
  \norm{A D_N^{1/2}\rho(t)} \leq s, \,\, \textrm{ and }\,\, \norm{A D_N^{1/2}\bw(t)}
  \leq s
\end{equation*}
for every $t\geq t_s$.

\section{Existence of an inertial manifold} \label{InertialMan}
%
To perform our analysis, we find convenient to rewrite system \eqref{eq.w}--\eqref{eq.div} as
\begin{align*}
&\partial_t \bar \bv - \nu \Delta \bar \bv + \nabla \cdot \fil{D_N \bar \bv 
	\otimes D_N \bar \bv} +\nabla \bar \pi = \bar \vartheta  \be_2\, , \\
&\partial_t \bar \vartheta  - k \Delta\bar \vartheta  + 
\nabla \cdot \fil{D_N\bar \vartheta  D_N\bar \bv}
= g\, ,  \\
&\nabla\cdot\bar \bv = 0\, \\
&\bw = A^{-1} \bv =: \overline{\bv}\\
&\rho = A^{-1} \vartheta =: \overline{\vartheta}
\end{align*}
and applying the operator $A$, term by term, to these equations we
get the equivalent problem
\begin{align}
&\partial_t  \bv - \nu \Delta  \bv + \nabla \cdot (D_N  \bar \bv 
\otimes D_N  \bar \bv) +\nabla  \pi =  \vartheta  \be_2\, , \label{eq.v}\\
&\partial_t  \vartheta  - k \Delta \vartheta  + \nabla \cdot (D_N \bar \vartheta
D_N \bar \bv)
= f\, , \label{eq.theta} \\
&\nabla\cdot \bv = 0,\, \label{eq.div-v}
\end{align}
with initial conditions
 \begin{equation} 
{\bv}(0, x) = {\bu}_0(x),\,\, 
\textrm{ and }\,\, {\vartheta}(0, t) =
  {\theta}_0(x), \label{data-1}
\end{equation}
and $f(\boldsymbol{x})=Ag(\boldsymbol{x})$ as introduced in \eqref{B-vero-2}.

Set $B(\bv, \bw):= \nabla \cdot (\bar \bv 
\otimes \bar \bw) $, $\mathcal R_1(\bv,\vartheta):=B(D_N \bar \bv,
D_N\bar \bv)$, 
$\mathcal R_2(\bv,\vartheta):=\nabla \cdot (D_N \bar \vartheta
D_N \bar \bv)$, $\mathcal R:=(\mathcal R_1,\mathcal R_2)$,
$\bV:=(\bv,\vartheta)\in (\cs L^2(\Omega))^3$, and 
$\eta=(\nu,\nu, k)$. 
Then, equations \eqref{eq.v}--\eqref{eq.theta} can be written as
\begin{equation*}
\partial_t \bV -\eta\cdot\Delta\bV+\mathcal R(\bV)=(\vartheta\ve e_2, f).
\end{equation*}

In order to prove the existence of the inertial manifold, we
show that $\mathcal R : (\cs L^2(\Omega))^3 \to (\cs L^2(\Omega))^3$ 
is locally Lipschitz continuous in $(\cs L^2(\Omega))^3$. 
We have that
\begin{equation} \label{Lip1}
\begin{aligned}
\|\mathcal{R}_1(\bV_1) - \mathcal{R}_1(\bV_2)\| = &\| B(D_N \bar \bv_1, D_N
\bar \bv_1) -  B(D_N \bar \bv_2, D_N \bar \bv_2)\|\\
\leq & \| B(D_N \bar \bv_1, D_N(
\bar \bv_1 - \bar \bv_2))\| + \| B(D_N (\bar \bv_1 - \bar \bv_2), D_N
\bar \bv_2 )\|\\
\leq& \|D_N \bar \bv_1\|_{L^4}\|\nabla D_N(
\bar \bv_1 - \bar \bv_2)\|_{L^4} + \| D_N(
\bar \bv_1 - \bar \bv_2)\|_{L^4}  \|\nabla D_N \bar \bv_2\|_{L^4}\\
\leq& c\lambda_1^{-1}(\|A D_N \bar \bv_1\| + \|A D_N \bar
\bv_2\|) \|A D_N( \bar \bv_1 - \bar \bv_2)\|\\
\leq& c\lambda_1^{-1}(\| D_N  \bv_1\| + \| D_N 
\bv_2\|) \| D_N(  \bv_1 -  \bv_2)\|\\
\leq& c\lambda_1^{-1}(N+1)^{3/2}(\| D_N^{1/2}  \bv_1\| + \| D_N^{1/2} 
\bv_2\|) \| \bV_1 -  \bV_2\|,
\end{aligned}
\end{equation}
and
\begin{equation} \label{Lip2}
\begin{aligned}
\|{\mathcal{R}}_2(\bV_1) -
{\mathcal{R}}_2(\bV_2)\| = &\| \nabla \cdot (D_N
\bar \vartheta_1 D_N \bar \bv_1) -  \nabla\cdot (D_N \bar \vartheta_2 D_N \bar \bv_2)\|\\
\leq & \| \nabla \cdot D_N (\bar \vartheta_1 - \bar \vartheta_2)D_N
\bar \bv_1 \| + \| \nabla \cdot D_N [\bar \vartheta_2  D_N
(\bar \bv_1 -\bar \bv_2 )]\|\\
\leq& c\lambda_1^{-1}\big(\|A  D_N \bar \bv_1\| \|A D_N( \bar
\vartheta_1 - \bar \vartheta_2)\|\\
&+ 
\| D_N  \vartheta_2 \| \| D_N(  \bv_1 - \bv_2)\|\big)\\
\leq & c\lambda_1^{-1}\big(\|  D_N  \bv_1\| \| D_N( 
\vartheta_1 - \vartheta_2)\|
+ \| D_N \vartheta_2 \| \| D_N( \bv_1 - \bv_2)\|\big)\\
\leq & c\lambda_1^{-1}(N+1)^{3/2}\big(\|  D_N^{1/2}  \bv_1\| \| 
\vartheta_1 - \vartheta_2\|\\
&+ \| D_N^{1/2} \vartheta_2 \| \| \bv_1 - \bv_2\|\big)\\
\leq &  c\lambda_1^{-1}(N+1)^{3/2}\big(\|  D_N^{1/2}  \bv_1\|
+ \| D_N^{1/2} \vartheta_2 \| \big) \| \bV_1 - \bV_2\|.
\end{aligned} \hspace{-1 cm}
\end{equation}

Thus, we have that
\begin{equation} \label{almost-global}
\|{\mathcal{R}}(\bV_1) -
{\mathcal{R}}(\bV_2)\|\leq 
c\lambda_1^{-1}(N+1)^{3/2}\big(2\|  D_N^{1/2}  \bv_1\| + \|  D_N^{1/2}  \bv_2\| 
+ \| D_N^{1/2} \vartheta_2 \| \big) \| \bV_1 - \bV_2\|,
\end{equation}
and we set $\mathfrak{L}:=c\lambda_1^{-1}(N+1)^{3/2}$.

In reproducing the argument used in \cite{HGTiti} to prove point (i) in
  \cite[Proposition 3]{HGTiti} (see also \cite[Proposition 9]{HGTiti}), i.e., the ``\emph{cone invariance
  property}'', which is stated in Proposition~\ref{prop.gap}--(i) below,
we introduce the following truncated nonlinearities: 
\begin{equation} \label{truncated-nonl}
\begin{aligned}
\mathcal{F}_1(\bV) :=
\chi_{\varrho}(\norm{D_N^{1/2}\bV})[\mathcal{R}_1(& \bV)
-\vartheta e_2],\,\, \mathcal{F}_2(\bV) := 
\chi_{\varrho}(\norm{D_N^{1/2}\bV})[\mathcal{R}_2(\bV)-f], \\[0.8 em]
&\textrm{ and }\, \,\, \mathcal F := (\mathcal F_1, \mathcal F_2).
\end{aligned}
\end{equation}
Here, $\chi$ is a smooth cut-off function
outside the ball of radius $\rho=2$ in $(\cs L^2(\Omega))^3$. Indeed,
let $\chi : \mathbb{R}^+ \to [0,1]$ with $\chi(r)=1$ for $0\leq r \leq 1$,
$\chi(r)=0$ for $r\geq 2$, and $|\chi'(r)|\leq 2$ for $r\geq 0$ (which
is null outside the ball of radius $\rho=2\tilde\varrho$).  
Define $\chi_{\tilde \varrho}(r):=\chi(r/\tilde \varrho)$ for
$r\geq0$.

Consider
\begin{align}
&\partial_t  \bv - \nu \Delta  \bv + 
  \chi_{\tilde{\varrho}}(\|D_N^{1/2}\bV\|)[\mathcal{R}_1(\bV)-\vartheta  \be_2]
  =  0\, , \label{eq.v1}\\
  &\partial_t  \vartheta  - k \Delta \vartheta  +
    \chi_{\tilde{\varrho}}(\|D_N^{1/2}\bV\|)[{\mathcal{R}}_2(\bV)
- f] = 0\, , \label{eq.theta1}
\end{align}
or equivalently
\begin{equation} \label{compact-form}
\partial_t \bV +\eta \cdot A\bV+\mathcal F(\bV)=0
\end{equation}
and observe that $\|\bv\|\cong \|A^{-1/2}D_N^{1/2}\bv\|\cong
\|A^{1/2}D_N^{1/2}\bar \bv\|\cong \|A^{1/2}D_N^{1/2}\bw\|$ and
  similarly we have that $\|\vartheta\|\cong
  \|A^{-1/2}D_N^{1/2}\vartheta\|\cong
\|A^{1/2}D_N^{1/2}\bar \vartheta\|\cong \|A^{1/2}D_N^{1/2}\rho\|$.

Let us recall that in the Section~\ref{estimates} (see the formulas in
  \eqref{eq:time-r}), we have shown that $\|D_N^{1/2}\bV (t)
\|\leq C\|A^{1/2}D_N^{1/2}\bV (t) \|\leq r=:\tilde{\varrho}$ for sufficiently large time $t \geq
t_{\tilde{\varrho}}$, where $r=\tilde \varrho$ is a suitable fixed
radius and $C=C(\alpha, \lambda_1)$.
Now, since $\mathcal{R} : (\cs L^2(\Omega))^3 \to   (\cs L^2(\Omega))^3$ is
locally Lipschitz, with calculations analogous to those
in \eqref{Lip1}--\eqref{almost-global},  it can be proved that
the truncated nonlinearity $\mathcal{F}$, defined in
\eqref{truncated-nonl}, is globally Lipschitz continuous
with Lipschitz constant given by $\mathcal{L}:=
c\lambda_1^{-1}(N+1)^{3/2}\tilde{\varrho}$.

\subsection{Spectral gap condition and inertial
  manifolds} \label{subsec.spectral} 
We now give the elements to prove our main result, i.e. the existence
of an inertial manifold to the system \eqref{eq.v}--\eqref{data-1}
(which is equivalent to \eqref{eq.w}--\eqref{eq.div}) in the sense of
Theorem~\ref{thm-utility} below.

Since the nonlinearity of \eqref{eq.v1}--\eqref{eq.theta1} is globally
Lipschitz, we will prove that this system has the
``\emph{strong squeezing property'}', stated precisely in Theorem~\ref{thm-utility},
provided that an appropriate ``\emph{spectral gap condition}'' (see
Proposition~\ref{prop.gap}--(i), below) is verified. 
Indeed, for $\gamma > 0$ and $n \in \N$, we introduce the cone
 \begin{equation} \label{cone}
\begin{aligned}
\Gamma_{n,\gamma} := \bigg\{ \begin{pmatrix} \bV_1\\ \bV_2
  \end{pmatrix}= \begin{pmatrix} (\bv_1, \vartheta_1)\\ (\bv_2,
    \vartheta_2) \end{pmatrix}  \in
  \mathcal{H}\X\mathcal{H} \, \colon \, \| D_N^{1/2} & Q_n\big((\bv_1, \vartheta_1) -(\bv_2,
  \vartheta_2)\big)\|\\
&\leq \gamma
\| D_N^{1/2}  P_n\big((\bv_1, \vartheta_1) -(\bv_2,
  \vartheta_2)\big)\|\bigg\} .
\end{aligned}
\end{equation} 
Here, we denote by $P_n$ the orthogonal projection from
$\mathcal{H}=H_0\X H^0$ onto $\textrm{span}\{(\phi_i, \varphi_j)\}$, 
where $\{\phi_i\}_{i=1}^\infty$ and $\{\varphi_j\}_{j=1}^\infty$ are
 orthonormal bases of $H_0$ and $H^0$, respectively, and also
define $Q_n := I - P_n$. 

Loosely speaking, the strong squeezing property asserts that if the
dynamics of two trajectories starts inside the cone $\Gamma_{n,\gamma}$, then the
trajectories stay inside the cone forever, and the higher Fourier
modes of the difference are dominated by the lower modes (i.e., the
cone invariance property); on the other hand, for as long as the two
trajectories are outside the cone, then the higher Fourier modes of
the difference decay exponentially fast, i.e., the decay
property (see Theorem~\ref{thm-utility} below for further details). 

Accordingly, in the case of the system \eqref{eq.v1}--\eqref{eq.theta1}, we have the following result.

\begin{proposition} \label{prop.gap} Let $\mathcal{L}
= c\lambda_1^{-1}(N+1)^{3/2}\tilde{\varrho}$ be the Lipschitz constant
associated with the truncated nonlinearity $\mathcal{F}$ in \eqref{truncated-nonl}.
Let $\bV_1=({\ve v}_1, \vartheta_1)$ and $\bV_2=({\ve v}_2, \vartheta_2)$ be two
solutions of \eqref{eq.v1}--\eqref{eq.theta1}. Then this system satisfies the
following properties.
\begin{itemize}
\item[\textrm{(i)}]
\textbf{The cone invariance property:} Assume that n is large enough such that
the ``spectral gap condition'' $\lambda_{n+1} - \lambda_n >
\frac{\mathcal{L}(\gamma+1)^2}{\nu\gamma}$ 
holds true. If  $\bv_1(t_0)  \in \Gamma_{n,\gamma}$ for some $t_0 \geq 0$,
then we have 
$\displaystyle \begin{pmatrix} (\bv_1(t), \vartheta_1(t))\\ (\bv_2(t), \vartheta_2(t)) \end{pmatrix}\in
\Gamma_{n,\gamma}$  for all $t\geq 0$.
\item[\textrm{(ii)}] \textbf{The decay property:} Assume that $n$ is large
  enough such that 
$\lambda_{n+1} >\nu^{-1}\mathcal{L}(\frac{1}{\gamma} + 1)$.
 If $\displaystyle \begin{pmatrix} (\bv_1(t), \vartheta_1(t))\\ (\bv_2(t), \vartheta_2(t)) \end{pmatrix}
\in \Gamma_{n,\gamma}$ for $0\leq t\leq T$, then it holds 
 \begin{equation*} \begin{aligned}
\big\| D_N^{1/2}  Q_n \big ((\bv_1(t), \vartheta_1(t))- (\bv_2(t),
   \vartheta_2(t))\big)\big\|
\leq
\big\| D_N^{1/2}  Q_n\big((\bv_1(0), \vartheta_1(0))- (\bv_2(0),
\vartheta_2(0))\big)\big\|e^{-\beta_n t}
\end{aligned}
\end{equation*}
for $0\leq t\leq T$, 
\end{itemize}
where $\beta_n =  \lambda_{n+1} |\eta|  - \mathcal{L} \big( \frac{1}{\gamma} +1\big)$
\end{proposition}
The proof of this proposition will be provided later in
Subsection~\ref{ssec.proof}.

Now, for the case at hand, we give the precise notion of Inertial manifold.

\begin{definition}[Inertial manifold] Consider the solution operator
 $S(t)$ associated with system \eqref{eq.v1}--\eqref{eq.theta1}. 
A subset $\mathcal{M} \subseteq \mathcal{H}$ is
 called an initial manifold for \eqref{eq.v1}--\eqref{eq.theta1} if the following properties are satisfied:
\begin{itemize}
\item[\textrm{(i)}] $\mathcal M$ is a finite-dimensional Lipschitz manifold;
\item[\textrm{(ii)}] $\mathcal M$ is invariant, i.e., $S(t)\mathcal{M}
  \subset
 \mathcal{M}$ for all $t \leq 0$;
\item[\textrm{(iii)}] $\mathcal M$ attracts exponentially all the
  solutions 
of \eqref{eq.v1}--\eqref{eq.theta1}; i.e.,
\begin{equation} \label{decay} \textrm{dist}(S(t)(\bv_0,\vartheta_0),
  \mathcal M) \to 0\,
  \textrm{ as }\, t \to \infty 
\end{equation}
 for every $(\bv_0, \vartheta_0) \in \mathcal{H}$, and the rate of decay in
 \eqref{decay} is exponential, 
uniformly for $(\bv_0, \vartheta_0)$ in bounded sets in $\mathcal{H}$.
\end{itemize}
\end{definition}
Property (iii) implies that $\mathcal M$ contains the global
attractor.
\smallskip

Next, we state a fundamental theorem concerning the fact that the strong 
squeezing property implies the existence of an inertial manifold 
 for dissipative evolution
equations. There are a number of proofs of such a result that can be
found, for instance, 
in \cite{Constantin, F-4, F-3}.  

\begin{theorem} \label{thm-utility}
Consider a nonlinear evolutionary equation of the type $v_t + Av + N
(v) = 0$, where $A$ is a linear, unbounded self-adjoint positive
operator, acting in a Hilbert space $H$, such that $A^{-1}$ is
compact and $N : H \to H$ is a nonlinear operator. Assume that the
solution $v(t)$ enters a ball in $H$ with radius $\varrho$ for
sufficiently large time $t$. For $\gamma > 0$ and $n \in\N$, 
we define the cone 
\begin{equation} \label{cone-H}
\Gamma_{n,\gamma} := \left\{ \begin{pmatrix} v_1\\ v_2\end{pmatrix}  \in
 {H}\X {H} \, \colon \, \big|Q_n(v_1 -v_2\big)\big|_H\leq \gamma
  \big|P_n(v_1 -v_2)\big|_H\right\} .
\end{equation}
Let $\chi$ be a smooth cut-off function
outside the ball of radius $2\varrho$ in $H$ and
let $\chi_{\varrho}$ be defined as in \eqref{eq.v1}--\eqref{eq.theta1}.
Assume that there exists $n \in \N$ such that the ``prepared'' equation 
\begin{equation} \label{prepared}
v_t + Av + \chi_{\varrho}(|v|)N (v) = 0
\end{equation}
satisfies the ``strong squeezing
property'', i.e., for any two solutions $v_1$ and $v_2$ of the
``prepared'' equation \eqref{prepared}, it holds that
\begin{itemize}
\item[(I)] if  $\begin{pmatrix} v_1(t_0)\\ v_2(t_0)\end{pmatrix}  \in
  \Gamma_{n,\gamma}$ for some $t_0 \geq 0$,
  then  $\begin{pmatrix} v_1(t)\\ v_2(t)\end{pmatrix} \in
  \Gamma_{n,\gamma}$ for all $t \geq t_0$;
\item[(II)] if  $\begin{pmatrix} v_1(t)\\ v_2(t)\end{pmatrix} \notin
  \Gamma_{n,\gamma}$ for $0\leq t\leq T$,
 then there exists $a_n >0$ such that
\begin{equation*}
|Q_n(v_1(t) - v_2(t))|_H \leq |Q_n(v_1(0)) - v_2(0)|_H e^{a_nt},\,\,
\textrm{ for }\,\, 0 \leq t \leq T.
\end{equation*}
\end{itemize}
Then, the ``prepared'' equation \eqref{prepared} admits an $n$-dimensional inertial
manifold in $H$. In addition, the following exponential tracking
property holds: 
for any $v_0 \in H$, there exists a time $\tau \geq 0$ and a solution 
$S(t)\varphi_0$ on the inertial manifold such that
\begin{equation*}
|S(t + \tau)v_0 - S(t)\varphi_0|_H \leq Ce^{-a_nt}, 
\end{equation*}
where the constant $C$ depends on $|S(\tau)v_0|_H$ and $|\varphi_0|_H$.
\end{theorem}

As a consequence of Proposition~\ref{prop.gap},  the strong squeezing property is verified 
 provided that $n\in\N$ is large enough; then in view of Theorem~\ref{thm-utility}, we have the
 following result

\begin{theorem} \label{thm.main} System \eqref{eq.v1}--\eqref{eq.theta1} admits a finite
  dimensional inertial manifold $\mathcal{M}$ in $\mathcal{H}$, i.e.,
  the solution $S(t)(\bu_0, \theta_0)$ of   \eqref{eq.v1}--\eqref{eq.theta1}
approaches the invariant Lipschitz manifold $\mathcal{M}$
exponentially fast in time. 
\end{theorem}

\subsection{Proof of Proposition~\ref{prop.gap}} \label{ssec.proof}
Let $\bV_1=(\bv_1, \vartheta_1)$ and $\bV_2=(\bv_2, \vartheta_2)$ be  two solutions
of \eqref{eq.v1}--\eqref{eq.theta1}. To show the cone invariance, i.e.
the condition given in Proposition~\ref{prop.gap}--(i) 
(see also \cite[Proposition~3-(i)]{HGTiti}), it is enough to show that  
$\begin{pmatrix} \bV_1(t)\\ \bV_2(t)\end{pmatrix}  =\begin{pmatrix}
  (\bv_1(t), \vartheta_1(t)) \\  
 (\bv_2(t), \vartheta_2(t)) \end{pmatrix}$
 cannot pass through the boundary of the cone if the dynamics starts
 inside the cone. More precisely, we will prove that
\begin{equation} \label{diff-derivative}
 \frac{d}{dt} \Big(\|D_N^{1/2}  Q_n\big(\bV_1(t) - \bV_2(t)\big)\| - \gamma
 \|D_N^{1/2}  P_n\big(\bV_1(t)  -
\bV_2(t))\|\Big) < 0 
\end{equation}
whenever  $\begin{pmatrix} \bV_1(t)\\ \bV_2(t)\end{pmatrix}  
\in \D \tilde{\Gamma}_{n, \gamma}$, where
$ \D \tilde{\Gamma}_{n, \gamma}$ stands for the boundary of the cone  
$\tilde{\Gamma}_{n, \gamma}$.

From \eqref{eq.v1}--\eqref{eq.theta1}, taking the difference
equations of $\bV_1$ and $\bV_2$,  we have that
\begin{align*}
&\partial_t(\bv_1 -\bv_2)-\nu \Delta(\bv_1 -\bv_2)+\mathcal{F}_1(\bV_1)-\mathcal{F}_1(\bV_2) = 0,\\
&\partial_t  (\vartheta_1 -\vartheta_2)  - k \Delta (\vartheta_1 -\vartheta_2)  +
  \mathcal{F}_2(\bV_1)-\mathcal{F}_2(\bV_2)   = 0\, , 
\end{align*}
or more compactly
\begin{equation} \label{compact-2}
\partial_t (\bV_1 - \bV_2) +\eta \cdot A(\bV_1-\bV_2)+\mathcal{F}(\bV_1) - \mathcal{F}(\bV_2)=0.
\end{equation}

By setting $\bp := P_n(\bV_1 -\bV_2)$, and
$\bq := Q_n(\bV_1 -\bV_2)$, 
we obtain
\begin{align}
& \bp_t +\eta\cdot A\bp+P_n\big(\mathcal{F}( \bV_1)-\mathcal{F}(\bV_2)\big) =0,
 \label{eq:A1}\\
& \bq_t +\eta\cdot A\bq+Q_n\big(\mathcal{F}(\bV_1)-
  \mathcal{F}(\bV_2)\big)= 0. \label{eq.A2}
\end{align}
Take the $L^2$-product of \eqref{eq:A1} with $D_N \bp$, to get
\begin{equation}
\frac{1}{2} \frac{d}{dt} \|D_N^{1/2} \bp\|^2 +
|\eta| \|A^{1/2}D_N^{1/2}\bp\|^2 + \big( D_N^{1/2} P_N
(\mathcal{F}(\bV_1) - \mathcal{F}(\bV_2)),
D_N^{1/2}\bp\big) 
=0\!\!\!\!
\end{equation}
Thus by \eqref{almost-global}, which is now a global property due to the
presence of $\chi_{\varrho}$, we have (with $\lambda_n$ 
large enough), we find
\begin{equation} \label{adattamento}
\begin{aligned}
\frac{1}{2} \frac{d}{dt} \|D_N^{1/2} \bp\|^2 \geq &
-\lambda_n |\eta|\|D_N^{1/2}\bp\|^2 - c\|D_N^{1/2}\big(\mathcal{F}(\bV_1) - \mathcal{F}(\bV_2)\big)\|
\|D_N^{1/2}\bp\| \\
\geq &
-\lambda_n |\eta|\|D_N^{1/2}\bp\|^2 - \tilde c \| D_N^{1/2}\big(\bV_1 -  \bV_2\big)\|
\|D_N^{1/2}\bp\| 
\end{aligned} 
\end{equation}
with $\tilde c=
c\lambda_1^{-1}(N+1)^{3/2}\tilde{\varrho}=\mathcal{L}$. 

Without loss of generality, we can assume $\|D_N^{1/2}\bp(t)\|>0$.
Otherwise, if $\|D_N^{1/2}\bp(t^\ast)\| = 0$ for some $t^\ast$, 
 then since we consider the boundary of the cone, we can assume 
$\|D_N^{1/2}\bq(t^\ast)\| = \gamma \|D_N^{1/2}\bp(t^\ast)\| = 0$, and thus
$D_N^{1/2}\bV_1(t^\ast) = D_N^{1/2}\bV_2(t^\ast)$. By the uniqueness of solutions, we
obtain $D_N^{1/2}\bV_1(t) = D_N^{1/2}\bV_2(t)$, for all $t \geq t^\ast$, and
the cone invariance property follows. Thus  $\|D_N^{1/2}\bp(t)\| > 0$, 
and we can divide both sides of \eqref{adattamento} by
$\|D_N^{1/2} \bp(t)\|$,
to reach
\begin{equation} \label{eq.bp}
 \frac{d}{dt} \|D_N^{1/2} \bp\| \geq - 
\lambda_n|\eta| \|D_N^{1/2}\bp\| - \mathcal{L} \|D_N^{1/2}\big( \bV_1 -  \bV_2\big)\|
\end{equation}
Similarly, taking the $L^2$-inner product of \eqref{eq.A2} 
against $D_N\bq$, we infer
\begin{equation} \label{eq.bq}
 \frac{d}{dt} \|D_N^{1/2} \bq\| \leq - \lambda_{n+1}|\eta|
 \|D_N^{1/2}\bq\| + \mathcal{L} \| D_N^{1/2}\big(\bV_1 -  \bV_2\big)\|
\end{equation} 
Now, multiplying equation \eqref{eq.bp} by $\gamma$ and subtracting it from
\eqref{eq.bq}, we obtain, exploiting the fact that
$\bp + \bq = \bV_1 - \bV_2$, the following relation
\begin{equation}
\begin{aligned}
\frac{d}{dt}\big( \|D_N^{1/2} \bq\|  -\gamma
\|D_N^{1/2} \bp\|\big) 
\leq |\eta| \big(\lambda_n
\gamma  \|D_N^{1/2} \bp\| - &\lambda_{n+1} \|D_N^{1/2} \bq\|
\big) \\
+ \mathcal{L} (\gamma + 1)(\| & D_N^{1/2} \bp\| + \|D_N^{1/2} \bq\|).
\end{aligned}
\end{equation}
Thus,  whenever $\|D_N^{1/2}\bq(t)\| = \gamma \|D_N^{1/2}\bp(t)\|$, 
that is 
$\begin{pmatrix} \bV_1(t)\\ \bV_2(t)\end{pmatrix}  
\in \D \tilde{\Gamma}_{n, \gamma}$, we have that
\begin{equation*}
\frac{d}{dt}\big( \|D_N^{1/2} \bq\|  -\gamma
\|D_N^{1/2} \bp\|\big) \leq  \left(|\eta| \big(\lambda_n
 - \lambda_{n+1}\big) + \mathcal{L}\frac{(\gamma
   +1)^2}{\gamma}\right) \|D_N^{1/2} \bq\| <0
\end{equation*}
due to the assumption $\lambda_{n+1} - \lambda_n > {\mathcal{L}(\gamma+1)^2}/{|\eta|\gamma}$ .
Whence
\begin{equation*}
\frac{d}{dt}\big( \|D_N^{1/2} \bq\|  -\gamma
\|D_N^{1/2} \bp\|\big) <0. 
\end{equation*}
this concludes point (i).

Now, to prove the decay property (ii), assume that $\begin{pmatrix}
  \bV_1(t)\\ \bV_2(t) \end{pmatrix}  \notin \tilde\Gamma_{n,\gamma}$ for
$0 \leq t \leq T$, then $\|D_N^{1/2}\bq(t)\| > \gamma\|
D_N^{1/2} \bp(t)\|$ for $0 \leq t \leq T$, and from
\eqref{eq.bq} it follows that
\begin{equation}
\begin{aligned}
 \frac{d}{dt} \|D_N^{1/2} \bq\|& \leq - \lambda_{n+1}|\eta|
 \|D_N^{1/2}\bq\|
 + \tilde c \big(\| \bp\| + \| \bq\|\big)\\
 &\leq - \lambda_{n+1}|\eta|
 \|D_N^{1/2}\bq\|
 + \hat c \big(\|D_N^{1/2}  \bp\| + \| D_N^{1/2}
 \bq\|\big)\\
& \leq -\Big[  \lambda_{n+1} |\eta|  -\mathcal{L}\big( \frac{1}{\gamma} +1\big)\Big]
\| D_N^{1/2} q\|=-\beta_n\|D_N^{1/2} q\|, 
\end{aligned}
\end{equation}
for $0\leq t\leq T$, where $\beta_n =  \lambda_{n+1} |\eta|  -
\mathcal{L}
\big( \frac{1}{\gamma} +1\big)$.
By Gronwall's inequality, we have that 
\begin{equation*}
   \|D_N^{1/2} \bq(t)\| \leq \|D_N^{1/2}
   \bq(0)\| e^{-\beta_n t}, \,\,
\textrm{ for }\,\, 0 \leq t \leq T.
 \end{equation*}
This concludes the proof.
\medskip
 
\noindent \textbf{Acknowledgments.} The authors are members of the Grup\-po Na\-zio\-na\-le per l'Ana\-li\-si Ma\-te\-ma\-ti\-ca, 
la Pro\-ba\-bi\-li\-t\`a e le lo\-ro Ap\-pli\-ca\-zio\-ni (GNAMPA) of the I\-sti\-tu\-to 
Na\-zio\-na\-le di Al\-ta Ma\-te\-ma\-ti\-ca (INdAM).

\end{document}